\documentclass[11pt, a4paper]{article}
\usepackage[english]{babel}
\usepackage{array}
\usepackage{mathdots}
\usepackage{amssymb,amsmath,amscd,amsthm}
\usepackage{mathtools}
\usepackage{thmtools}
\usepackage{xcolor}
\usepackage{comment}
\usepackage{bold-extra} 
\usepackage{hyperref} 
\usepackage{bbm} 
\usepackage[title]{appendix} 
\usepackage[round]{natbib}

\declaretheoremstyle[
spaceabove=10pt, spacebelow=10pt,
headfont=\normalfont\scshape\bfseries,
notefont=\normalfont\scshape,
notebraces={(}{)},
bodyfont=\itshape,
headpunct=,
postheadspace=1em
]{THMS}

\declaretheoremstyle[
spaceabove=6pt, spacebelow=6pt,
headfont=\normalfont\scshape\bfseries,
notefont=\normalfont\scshape,
notebraces={(}{)},
bodyfont=\normalfont,
headpunct=,
postheadspace=1em
]{DEF}

\declaretheoremstyle[
spaceabove=12pt, spacebelow=6pt,
headfont=\normalfont\scshape\bfseries,
notefont=\normalfont\scshape,
notebraces={(}{)},
bodyfont=\normalfont,
headpunct=,
postheadspace=1em,
qed=$\blacksquare$
]{EX}

\declaretheoremstyle[
spaceabove=12pt, spacebelow=12pt,
headfont=\normalfont\scshape\bfseries,
notefont=\normalfont\scshape\bfseries,
notebraces={}{},
bodyfont=\normalfont,
headpunct={:},
qed=$\blacksquare$,
postheadspace=1em
]{PROOF}

\declaretheorem[style=THMS, numberwithin=section, name=Theorem]{thm}
\declaretheorem[style=THMS, numberwithin=section, name=Proposition]{prop}
\declaretheorem[style=THMS, numberwithin=section, name=Lemma]{lma}

\declaretheorem[style=DEF, numbered=no, name=Notation]{notation}

\declaretheorem[style=EX, name=Example]{ex}

\declaretheorem[style=PROOF, numbered=no, name=Proof of ]{proofc} 

\numberwithin{equation}{section}

\newcommand*{\E}{\operatorname{E}}
\newcommand{\norm}[1]{\left\lVert#1\right\rVert}

\mathtoolsset{showonlyrefs} 

\title{Semiparametric Fisher Information in Models parametrized by a Normed Space}
\author{Telmo P\'erez-Izquierdo\footnote{University of the Basque Country (UPV/EHU). Email address: telmo.perezizquierdo@ehu.eus.} }
\date{\today}

\begin{document}

\maketitle

\begin{abstract} 

This paper studies semiparametric Fisher information in models parametrized by general normed spaces. The main contribution is to establish that positive semiparametric Fisher information is equivalent to the gradient of the parameter of interest lying in the range of the adjoint score operator. This result generalizes a key theorem \citet{van1991differentiable} and provides a unified framework linking differentiability and information, beyond Hilbert spaces. The paper develops a normed-space mean-square-differentiable models for two canonical problems: estimation of the average of a known transformation and estimation of a density at a point. In these applications, it shows that positive information holds if and only if the transformation has finite variance and if and only if the density has positive mass at the evaluation point, respectively. These findings offer a novel information-theoretic perspective on known minimax results and clarify the conditions under which root-n estimation is possible.

\vspace{11pt}

\end{abstract}

\section{Introduction}

The problem of asymptotic estimation and inference can be framed in a general setup where the goal is to recover a parameter $\beta_0 \equiv \kappa(P_0)$, for a functional $\kappa\colon \mathcal{P} \to \mathbb{R}$ and a collection of distributions $\mathcal{P}$, where $P_0 \in \mathcal{P}$ is the true distribution of the data. The key question is how well one can estimate $\beta_0$. In this setting, local asymptotic minimax and convolution/representation theorems give lower bounds for the asymptotic variance of regular estimators of $\beta_0$ in locally asymptotically normal models \citep{le1972limits,hajek1970characterization,hajek1972local}. These results have been extended to nonparametric and semiparametric models \citep{koshevnik1976non,begun1983information, van1991differentiable}.

A key insight in estimation of $\beta_0$ is that of differentiability of the defining functional $\kappa$ \citep{pfanzagl2012contributions, van1991differentiable}. Differentiability is a necessary condition for regular estimation of $\beta_0$ at root-n rate \citep[Th.~2.1]{van1991differentiable}. Moreover, there is a stark link between root-n estimability, differentiability, and positive semiparametric Fisher information. For models $\mathcal{P} = \{P_\lambda\colon \lambda\in\Lambda\}$ parametrized by a subset $\Lambda$ of a Hilbert space, \citet{van1991differentiable} establishes that differentiability of $\kappa(P_\lambda)\equiv \psi(\lambda)$, for a functional $\psi\colon \Lambda \to \mathbb{R}$, is equivalent to positive semiparametric Fisher information for $\beta_0$. Hence, positive information is a necessary condition for root-n estimation of $\beta_0$ \citep[see also][]{chamberlain1986asymptotic}.   

The above framework is quite general since, given a dominating measure one can parametrize $\mathcal{P}$ by the square-roots of the densities, which live in a Hilbert space. However, key results rely on mean-square differentiability of the model, which requires a continuous score operator, and on pathwise differentiability of $\psi$, which requires a continuous gradient. These conditions tend to hold for well-behaved parameters, but may fail when one attempts to show that a parameter is not root-n estimable. In such cases, certain smoothness conditions on the model readily lead to mean-square differentiability and a pathwise differentiable functional. Nevertheless, the restricted model ends up being parametrized by a general normed space, not a Hilbert space.   

The contribution of this paper is to generalize the link between positive information and differentiability to models parametrized by normed spaces. I extend Theorem~4.1 in \citet{van1991differentiable} to normed spaces, showing that positive semiparametric Fisher information is equivalent to the gradient being in the range of the adjoint score operator. The usefulness of the extension is illustrated in two classical problems: estimation of the average of a known transformation and estimation of the density at a point. For both applications, I introduce a normed-space mean-square-differentiable model where the functionals of interest are pathwise differentiable with continuous gradients. I find that semiparametric Fisher information is positive for estimation of the average of a known transformation if and only if the transformation has finite variance. For estimation of the density at a point, semiparametric Fisher information is positive if and only if there is positive mass at the point of interest. Up to my knowledge, these results are novel, even if they are aligned with the literature. 

There is a large recent literature on estimation of the mean for heavy-tailed distributions, which is closely related to estimation of the average of a known transformation with  possibly infinite variance \citep[for an overview, see][]{lugosi2019mean}. For instance, \citet[Th.~3.1]{devroye2016sub} give a (non-asymptotic) minimax bound of $n^{-s/(1+s)}$ for estimation of the mean in the class of distributions with $\E_0[|g(X)|^{1+s}] \leq C$, with $s \in (0, 1]$. On the other hand, slower than root-n rates for density estimation are well established. For instance, \citet[Th.~5.1]{ibragimov1981problem} give a minimax bound for estimation of $n^{-s/(2s+1)}$, with $s \geq 1$ being the smoothness of the class of densities (p.~235). My contribution is to link these results with semiparametric information theory, providing a standard way to obtain them.

Zero information and non-differentiability results for the average of a known transformation and for the density at a point add up to the literature on efficient semiparametric estimation. \citet{chamberlain1986asymptotic} finds information bounds in semiparametric models with censoring. \citet{ritov1990achieving} discuss estimation of the squared integrated density, where information is positive but no root-n estimator exists. \citet{khan2010irregular} find zero information results for a semiparametric the binary choice model and treatment effects under exogenous selection. These authors use pathwise calculations to obtain their results. Closer to this paper, \citet{van1991differentiable} obtains results for mixing, censoring and truncation models. \citet{escanciano2018semiparametric} finds necessary and sufficient conditions for positive information for willigness to pay and average risk aversion. To obtain these results, they study the relative position of the gradient with respect to the range of the adjoint score operator.  

The rest of the paper is organized as follows. Section~\ref{sec:setup} introduces the setup and motivating examples. Section~\ref{sec:overview} gives an overview of the results. Section~\ref{sec:diff_info} contains the main result. Section~\ref{sec:applications} applies the main result to the examples. Section~\ref{sec:conclusion} concludes. All the proofs are in the Appendix.

\section{Setup and motivating examples}
\label{sec:setup}

Let $(V, \norm{\cdot})$ be a normed vector space. I consider a probability model over a measurable space $(\mathcal{X}, \mathcal{B})$ given by the collection of distributions $\mathcal{P} \equiv \{P_\lambda \colon \lambda \in \Lambda\}$, where $\Lambda \subseteq V$. The true distribution of the random element $X$ is $P_0 \in \mathcal{P}$. I assume that there is one parameter $\lambda_0 \in \Lambda$ that generated the data, i.e., $P_0=P_{\lambda_0}$ (well-specified model). Through this paper, it is assumed that the distributions in $\mathcal{P}$ are dominated by a $\sigma$-finite measure $\mu$. For each $\lambda\in\Lambda$, $p_\lambda$ denotes the density of $P_\lambda$ with respect to $\mu$ and $p_0\equiv p_{\lambda_0}$.

The setup here parallels that of \citet{van1991differentiable}. Consider the set of paths $t \mapsto \lambda_t$, $t \in (0, \varepsilon) \subset \mathbb{R}$ and $\varepsilon > 0$, satisfying the following regularity condition:
\begin{equation} \label{eq:reg_path}
    t^{-1} (\lambda_t - \lambda_0) \to \alpha \in V, \quad \text{as } t \downarrow 0.
\end{equation}
The collection of all the deviations $\alpha$ obtained from all the paths is the tangent space $\mathcal{T}$, which is a assumed to be a closed linear subspace of $V$. I also assume that the model is mean-square differentiable, that is, there exists a continuous linear operator $A\colon \mathcal{T} \to L_2(P_0)$, called the score operator, such that
\begin{equation} \label{eq:msd}
    \int \left( \frac{\sqrt{p_{\lambda_t}} - \sqrt{p_0}}{t} - \frac{1}{2} A\alpha \sqrt{p_0} \right)^2 d\mu \to 0, \quad \text{as } t \downarrow 0,
\end{equation}
for every path $t \mapsto \lambda_t$ satisfying \eqref{eq:reg_path}. For the measure $P_0$, the space $L_2(P_0)$ denotes the space of (equivalent classes of) square-integrable functions, i.e., functions $f\colon \mathcal{X} \to \mathbb{R}$ with $\norm{f}_2 < \infty$, where $\norm{f}_2 \equiv (\int f^2 dP_0)^{1/2}$. 

The goal is to estimate the parameter $\beta_0 \equiv \psi(\lambda_0)$, for a functional $\psi\colon\Lambda \to \mathbb{R}$. The functional $\psi$ is pathwise differentiable, that is, there exists a continuous linear operator $\dot{\psi}\colon \mathcal{T} \to \mathbb{R}$ such that
\begin{equation}
    t^{-1} \left( \psi(\lambda_t) - \psi(\lambda_0) \right) \to \dot{\psi}\alpha, \quad \text{as } t \downarrow 0,
\end{equation}
for every path $t \mapsto \lambda_t$ satisfying \eqref{eq:reg_path}. I refer to $\dot\psi$ as the gradient of $\psi$ at $\lambda_0$. 

I close this section with a couple of classical examples. They help to
illustrate the need to generalize from models parametrized by Hilbert spaces to cases where $V$ is just a normed space.

\begin{ex}[Average of a known transformation]
	\label{ex:mean}\label{ex_mean}

	Consider a measurable function $g\colon \mathcal{X} \to \mathbb{R}$, with $g \in L_q(P_0)$ for some $1 \leq q \leq 2$. The goal is to estimate $\beta_0 \equiv \E_0[g(X)]$, where the expectation is taken w.r.t. $P_0$, the true distribution of $X$. It is convenient to impose certain smoothness conditions on the distributions of the model, considering densities in $\mathcal{P} \equiv \{p \in L_1(\mu)\colon fp \in L_1(\mu), \forall f \in L_q(P_0)\}$. This requires that $\E_p[|f(X)|] = \int |f|pd\mu < \infty$ for every $f \in L_q(P_0)$, not just the transformtion $g$ of interest. The additional smoothness is enough to guarantee mean-square differentiability of the model (see Section~\ref{sec:applications}).  

	 I reparametrize the problem in terms of deviations. Suppose that $\lambda \in L_1(\mu)$ measures a perturbation to $p_0$, in the sense that each $\lambda$ leads to a density $p_\lambda \equiv p_0\cdot (1+\lambda)$. I parametrize the problem by the set of deviations $\lambda$ that lead to valid densities:
	\begin{equation}
		\Lambda \equiv \{ \lambda \in L_1(\mu) \colon fp_0(1+\lambda) \in L_1(\mu), \forall f \in L_q(P_0) \}.
	\end{equation}
	It is convenient to express $\Lambda$ in the following equivalent form (see Section~\ref{sec:applications}):
	\begin{equation}
		\Lambda = \{ \lambda \in L_1(p_0) \colon \E_0[|f\lambda|] < \infty , \forall f \in L_q(P_0) \}.
	\end{equation}

	There are several advantages to working in this parameter space. First, note that the ``true parameter" is $\lambda_0=0$, that is, no deviation from the ``true distribution" $p_0$. More importantly, one  can characterize which functions live in $\Lambda$, depending on the smoothness of $g \in L_q(P_0)$. A standard duality argument shows that $\Lambda = L_q(P_0)^* = L_{q'}(P_0)$, with  $1/q+1/q'=1$ (see Proposition~\ref{prop:duality}). This highlights that, when $g$ is less smooth (smaller $q$), one must consider smoother (greater $q'$) deviations. Note that $\Lambda = L_2(P_0)$, a Hilbert space, only when $q=2$, that is, when $g$ has finite variance.

	In this example, I consider $V=\Lambda=L_{q'}(P_0)$. The tangent space is the whole parameter space: $\mathcal{T} = L_{q'}(P_0)$. Indeed, $L_{q'}(P_0)$ is complete and for any $\alpha \in L_{q'}(P_0)$, the path $\lambda_t = t\alpha$ satisfies \eqref{eq:reg_path}. The score operator is the inclusion map from $L_{q'}(P_0)$ to $L_2(P_0)$, that is, $A\alpha = \alpha$ (see Proposition~\ref{prop:est_mean_score}).  The parameter of interest is given by the functional $\psi(\lambda) \equiv \int g(1+\lambda)p_0d\mu$, with gradient $\dot{\psi}\colon L_{q'}(P_0) \to \mathbb{R}$ given by $\dot{\psi}\alpha = \int \alpha g p_0 d\mu$ (see Section~\ref{sec:applications}). Realize that the gradient is not $L_2(P_0)$-continuous when $q<2$. Here, the extension to normed spaces is crucial since continuity of $\dot\psi$ requires to restrict the tangent space to $L_{q'}(P_0)$, which is not a Hilbert space when $q<2$.
	
	Estimation of the average of a known transformation, with $q \in (1, 2]$, serves a running example to illustrate the results. The case of $g \in L_1(P_0)$ is presented in Section~\ref{sec:applications}.
\end{ex} 

\begin{ex}[Density at a point]
	\label{ex:density}

	Consider that $\mathcal{X}$ is a locally compact Hausdorff space. Let $P_0$ be the distribution of a random element $X \in \mathcal{X}$, with pdf $p_0$ with respect to a $\sigma$-finite measure $\mu$. Assume further that $\mu$ is locally finite. The parameter of interest is $\beta_{0}=p_0(x)$ for a point $x\in\operatorname{supp}(P_0)$, where $\operatorname{supp}(P_0) \equiv \{x \in \mathcal{X} \colon p_0(x) > 0\}$ is the (interior of the) support of $P_0$. To avoid issues with equivalence classes of functions, the probability model $\mathcal{P}$ is restricted to measures with continuous densities with respect to~$\mu$: $\mathcal{P} \equiv \mathcal{C}(\mathcal{X})$, where $\mathcal{C}(\mathcal{X})$ is the set of continuous functions $p \colon \mathcal{X} \to \mathbb{R}$.

	I reparametrize the problem in terms of local deviations. The construction is based on a compact neighborhood $K$ of $x$, since the behavior of densities outside this neighborhood is irrelevant for estimation of $p_0(x)$. The parameter of the model is a deviation in $\Lambda \equiv \mathcal{C}(K)$, which translates into a density $p_\lambda = p_0 + u \lambda$, where $u$ is a continuous function equal to one in a neighborhood of $x$ and vanishing outside $K$. Note that $\lambda_0 = 0$. Even if this parametrization does not include all the densities in $\mathcal{P}$, it is rich enough to capture their local behavior. That is, for every $p\in\mathcal{P}$, there exists a $\lambda \in \Lambda$ such that $p_\lambda = p$ in a neighborhood of $x$ (see Section~\ref{sec:applications} for more details).
	
	I consider $V = \Lambda = \mathcal{C}(K)$ endowed with the supremum norm $\norm{\lambda}_\infty \equiv \sup_{x\in K} |\lambda(x)|$. As in Example~\ref{ex:mean}, since $\mathcal{C}(K)$ is complete, the tangent space is $\mathcal{T} = \mathcal{C}(K)$. The score operator $A\colon \mathcal{C}(K) \to L_2(P_0)$ is $A\alpha = u\alpha/p_0$. The advantage of the parametrization in terms of local deviations is that it guarantees mean-square differentiability of the model (see Proposition~\ref{prop:est_dens_score}). The parameter of interest is given by the functional $\psi(\lambda) = p_0(x) + u(x)\lambda(x)$, with gradient $\dot{\psi}\alpha = \alpha(x)$ (see Section~\ref{sec:applications}). Evaluation functionals are not $L_2(P_0)$-continuous, so the extension to normed spaces is crucial to analyze the problem.  
\end{ex}

\begin{notation}
I use the following notation throughout the paper. If $V$ and $W$ are normed spaces and $A\colon V \to W$ is a linear operator, $\mathcal{R}(A) \equiv \{ w \in W \colon \exists v \in V \text{ s.t. } Av = w \}$ denotes the range of $A$ and $\mathcal{N}(A) \equiv \{ v \in V \colon Av = 0 \}$ denotes the kernel of $A$. If $E$ is another normed space and $T \colon W \to E$ is a linear operator, the composition of $A$ and $T$ is $TA \colon V \to E$, defined by $[TA]v = T(Av)$ for every $v \in V$. It is often convenient (see below) to denote ``operator $A$ evaluated at point $v$" by $\langle v, A \rangle$ instead of $Av$. For instance, $\langle v, TA \rangle \equiv  \langle Av, T \rangle$ defines the composition of $T$ and $A$.

For a normed space $V$, $V^*$ denotes its dual, that is, the space of continuous linear functionals on $V$. For a continuous linear operator $A\colon V \to W$, its adjoint $A^* \colon W^* \to V^*$ is defined by $\langle v, A^* w^* \rangle = \langle Av, w^* \rangle$ for every $v \in V$ and $w^* \in W^*$. This notation resembles that of Hilbert spaces, which are self-dual and where $\langle \cdot, \cdot \rangle$ can be understood as the inner product.
\end{notation}

\section{Overview of the results}
\label{sec:overview}

The main result provides a link between positive semiparametric information and the position of the gradient $\dot{\psi}$ in the dual of the tangent space $\mathcal{T}$. Semiparametric Fisher information for $\beta_0$ is defined as
\begin{equation}
    \mathcal{I} \equiv \inf_{\alpha \in \mathcal{T}} \mathcal{I}(\alpha), \text{ where} \quad \mathcal{I}(\alpha) \equiv \frac{\norm{A\alpha}_2^2}{|\dot{\psi}\alpha|^2}.
\end{equation}
$\mathcal{I}(\alpha)$ gives the information for $\beta_0$ in the one-dimensional parametric submodel $\{ p_{\lambda_t} \colon t \in (0, \varepsilon) \}$ generated by the a path $\lambda_t = \lambda_0 + t\alpha + o(t)$. $A\alpha$ is the score for $t$ in the submodel and $\dot{\psi}\alpha$ is the derivative of $\psi(\lambda_t)$ w.r.t. $t$. That is, $\mathcal{I}(\alpha)^{-1}$ is the usual Cram\'er-Rao bound for estimation of $\beta_0$ in the submodel \citep{cramer1999mathematical,rao1945information}. Semiparametric Fisher information $\mathcal{I}$ considers the worst-case scenario across all the submodels generated by paths satisfying \eqref{eq:reg_path}.  

I establish that $\mathcal{I}> 0$ if and only if $\dot{\psi} \in \mathcal{R}(A^*)$. This result extends the link between positive information and differentiability \citep[Th.~4.1]{van1991differentiable}, from models parametrized by Hilbert spaces to more general normed spaces. The extension is non-trivial, providing a proof of the ``only if'' part that does not rely on Hilbert-space techniques. The result is also valid for functionals $\psi\colon \Lambda \to W$, with $W$ an arbitrary normed linear space. In that case, one just requires that $w^*\dot\psi \in \mathcal{R}(A^*)$ for all $w^*  \dot\psi \colon \mathcal{T} \to \mathbb{R}$, with $w^*\in W^*$. That is, the condition becomes $\mathcal{R}(\dot\psi^*) \subseteq \mathcal{R}(A^*)$, for the adjoint $\dot\psi^* \colon W^* \to \mathcal{T}^*$.

The usefulness of the extension is illustrated in two classical problems. In Section~\ref{sec:setup}, I have introduced normed space models for estimation of the average of a known transformation and for estimation of the density at a point. These models are mean-square differentiable, as in equation~\eqref{eq:msd}. Moreover, the functional of interest is pathwise differentiable, as in equation~\eqref{eq:reg_path}. Pathwise differentiability of the functional is essential to apply the main result of the paper, as it is to apply Theorem~3.1 in \citet{van1991differentiable}. However, it does not generally hold for Hilbert parameters, so the extension to normed spaces is crucial to analyze these problems. 

I find that semiparametric Fisher information is positive for estimation of the average of a known transformation if and only if $g \in L_2(P_0)$, that is, the transformation has finite variance. For estimation of the density at a point, semiparametric Fisher information is positive if and only if $\mu(\{x\}) > 0$, that is, there is positive mass at the point of interest. The results are consistent with the literature. What is novel here is the modelling strategy, the link between positive information and differentiability, and the impossibility of root-n estimation when information is zero. 

\section{Differentiability and positive information}
\label{sec:diff_info}

The proof strategy followed here discloses the link between positive information and the fact that $\dot\psi\in\mathcal{R}(A^*)$. The strategy relies on two ideas. First, I show that positive information and $\dot\psi\in\mathcal{R}(A^*)$ are linked to smoothness of the problem. That is, both concepts capture how easy it is to recover $\beta_0=\psi(\lambda_0)$ from the data $p_\lambda$.

To see this, it is useful to simplify the problem and consider one-to-one score operators. This establishes a bijection between the tangent space of the data $\mathcal{R}(A)$ and the tangent space of the parameter $\mathcal{T}$. For each local deviation $\delta \in \mathcal{R}(A)$ in the data, there is a unique local deviation $\alpha$ in the parameter. Therefore, one can locally recover $\lambda_t = \lambda_0 + t\alpha + o(t)$ from the data, and the information in the model reduces to smoothness of the functionals involved.  

Of course, the one-to-one assumption on the score operator is very restrictive. The second idea on which the proof is based is the following: if a local identifiability condition holds \citep[see][]{van1991differentiable,escanciano2018semiparametric}, then $A$ is ``sufficiently invertible". By this I mean that $A$ is invertible enough to locally recover $\psi(\lambda_t)$. The precise meaning of this is made clear after the discussion of the one-to-one $A$ case.

\subsection{One-to-one score operators}

Since $A$ is one-to-one and onto $\mathcal{R}(A) \subseteq L_2(p_0)$, the inverse $A^{-1}\colon\mathcal{R}(A)\to \mathcal{T}$ is guaranteed to exists. It is trivial to check that $A^{-1}$ is linear. $A^{-1}$ translates local deviations in the data to local deviations in the parameter of the model. The linear functional $\dot\psi  A^{-1} \colon \mathcal{R}(A) \to \mathbb{R}$, in turn, translates local deviations in the data to local deviations in the parameter of interest. Positive information is linked to continuity of $\dot\psi A^{-1}$:
\begin{prop}{\label{prop:posinfo_continuous}}
	Let $A$ be one-to-one. Then,  $\mathcal{I}>0$ if and only if $\dot\psi A^{-1}$ is continuous.
\end{prop}

It is easier to grasp the above result in the simplified setting of $\mathcal{T} = L_2(p_0)$ and $A$ the identity. In this case, information in direction $\alpha$ reduces to
\begin{equation}
		\mathcal{I}(\alpha) =\frac{\norm{\alpha}^2_{2}}{|\dot{\psi}\alpha|^2} \geq \mathcal{I}.
\end{equation}
Then, positive information is equivalent to variation in $|\dot\psi\alpha|$ being bounded by $\mathcal{I}^{-1/2}\norm{\alpha}_{2}$, that is, continuity of the gradient $\dot\psi$.

The next proposition shows that the relative position of $\dot\psi$ with respect to the range of the adjoint score operator is also linked to continuity of $\dot\psi A^{-1}$:
\begin{prop}{\label{prop:range_continuous}}
		Let $A$ be one-to-one. Then,  $\dot\psi \in \mathcal{R}(A^*)$ if and only if $\dot\psi A^{-1}$ is continuous.
\end{prop}
This results can be combined with Proposition~\ref{prop:posinfo_continuous} to get that positive information is equivalent of $\dot\psi\in\mathcal{R}(A^*)$ for one-to-one score operators. The next section is devoted to sidestep the one-to-one assumption. Theorem~\ref{thm:info_range} states the main result.

\subsection{General score operators}

For general score operators, the implication $\dot\psi \in \mathcal{R}(A^*) \implies \mathcal{I}>0$ follows the argument in \citet[Th.~4.1]{van1991differentiable}. The reverse implication is more involved and it relies on collapsing all the paths that have the same score into equivalence classes. This is possible under the local identifiability condition \citep{van1991differentiable,escanciano2018semiparametric}:
\begin{equation}
	\mathcal{N}(A) \subseteq \mathcal{N}(\dot\psi).
\end{equation}

Since positive information implies local identifiability, the construction is valid to show that $\mathcal{I}>0$ implies $\dot\psi \in \mathcal{R}(A^*)$. Intuitively, if local identifiability fails, there is at least a path such that $\sqrt{p_{\lambda_t}} = \sqrt{p_0} + o(t)$ but $\psi(\lambda_t) = \beta_0 + t\dot\psi\alpha + o(t)$ with $\dot\psi\alpha \neq 0$. This is a path with zero information, since one cannot locally distinguish $p_{\lambda_t}$ from $p_0$ but $\psi(\lambda_t)$ is locally different from $\psi(\lambda_0)$.  

The construction is based on the quotient space $\mathcal{Q} \equiv \mathcal{T}/\mathcal{N}(A)$ containing all equivalence classes $[\alpha] \equiv \alpha+\mathcal{N}(A)=\{\alpha+n\colon n\in \mathcal{N}(A)\}$. Elements in $\mathcal{Q}$ are sets containing locally indistinguishable deviations: if $\alpha'\in [\alpha]$, then $A\alpha'=A\alpha$. The quotient space $\mathcal{Q}$ is, indeed, a normed vector space \citep[Th.~1.7.4]{megginson1998introduction}. There, one can define the quotient score operator $A_\mathcal{Q} \colon \mathcal{Q} \to L_2(p_0)$, given by $A_\mathcal{Q}[\alpha] \equiv A\alpha$. Note that $A_\mathcal{Q}$ is well defined by construction.

The proof relies on the following insights. First, the quotient score operator $A_\mathcal{Q}$ is one-to-one and onto $\mathcal{R}(A_\mathcal{Q}) = \mathcal{R}(A)$. Hence, the results in the preceeding section apply to $A_\mathcal{Q}$. Second, one can define the quotient gradient $\dot\psi_\mathcal{Q} \colon \mathcal{Q} \to \mathbb{R}$ as $\dot\psi_\mathcal{Q}[\alpha] \equiv \dot\psi\alpha$. Under the local indentifiability condition $\mathcal{N}(A) \subseteq \mathcal{N}(\dot\psi)$, the functional $\dot\psi_\mathcal{Q}$ is also well defined. Third, under the local identifiability condition, information along the submodel $\lambda_t = \lambda_0 + t\alpha + o(t)$ is the same as information along the submodel $\lambda_t' = \lambda_0 + t\alpha' + o(t)$, for every $\alpha'\in [\alpha]$. This means that, in essence, one can collect all the paths that have the same score into equivalence classes and work with the corresponding tangent space.
   
The next lemma formally states the above insights. 
\begin{lma}{\label{lma:properties}}
	Assume that $\mathcal{N}(A) \subseteq \mathcal{N}(\dot\psi)$. Then,
	\begin{enumerate}
		\item $A_\mathcal{Q}$ is well-defined, linear, continuous, one-to-one, and onto $\mathcal{R}(A)$. 
		
		\item $\dot\psi_\mathcal{Q}$ is well-defined, linear, and continuous.
		
		\item It holds that $\mathcal{I}=\mathcal{I}_\mathcal{Q}$, with 
		\begin{equation}
				\mathcal{I}_\mathcal{Q} \equiv \inf_{[\alpha]\in \mathcal{Q}} \frac{ \norm{A_\mathcal{Q}[\alpha]}_{2}^2}{\left\lvert\dot\psi_\mathcal{Q}[\alpha]\right\lvert^2}.
		\end{equation}
		
		\item It holds that $\dot\psi_\mathcal{Q} \in \mathcal{R}(A_\mathcal{Q}^*) \iff \dot\psi\in\mathcal{R}(A^*)$.
	\end{enumerate}
\end{lma}

	With these results, I am ready to state the main theorem of the paper. The proof strategy is to first show that positive information implies local identification. Then, I can apply Propositions~\ref{prop:posinfo_continuous} and \ref{prop:range_continuous} to $A_\mathcal{Q}$ and $\dot\psi_\mathcal{Q}$.  
	\begin{thm}{\label{thm:info_range}} $\mathcal{I}>0$ if and only if $\dot\psi\in \mathcal{R}(A^*)$.
	\end{thm}

	The next example shows how to apply the result of the theorem:
\begin{ex}[continues=ex:mean]
	One needs to study the adjoint score operator $A^* \colon L_2(P_0)^* \to L_{q'}(P_0)$. Since $L_2(P_0)$ is a Hilbert space, $L_2(P_0)^* = L_2(P_0)$, with $=$ indicating that there exists an isometric isomorphism. Here, I consider the case where $g \in L_q(P_0)$, with $q\in(1,2]$. Hence, $q' < \infty$, and the dual space is $L_{q'}(P_0)^* = L_q(P_0)$ \citep[Th.~6.16]{rudin1970real}. From the definition of the gradient $\dot\psi$, it follows that the representer of $\dot\psi$ on $L_q(P_0)$ is the function $g$. Since the score operator is the inclusion map, $\langle \alpha, A^*\delta \rangle = \langle A\alpha, \delta \rangle = \langle \alpha, \delta \rangle$ for every $\alpha \in L_{q'}(P_0)$ and $\delta \in L_{2}(P_0)$. From this it follows that $A^*\delta = \delta$, so that $A^*$ is the inclusion map from $L_2(P_0)$ to $L_{q}(P_0)$. Its range is $\mathcal{R}(A^*) = L_2(P_0) \subseteq L_q(P_0)$.

	When the known transformation has finite variance, i.e., $q=2$, the range of the adjoint score operator fills the whole space. Hence, $\dot\psi \in \mathcal{R}(A^*)$ and semiparametric information for $\beta_0 \equiv \E_0[g(X)]$ is positive according to Theorem~\ref{thm:info_range}. This is the classical Hilbert-space result. To obtain the negative result, that is, that information is negative without the finite variance assumption, one needs to extend the analysis to normed spaces. If $g \in L_{q}(P_0) \setminus L_2(P_0)$, then $\dot\psi \notin \mathcal{R}(A^*)$ and therefore, by Theorem~\ref{thm:info_range}, semiparametric information for $\beta_0$ is zero. In this case, Theorem~3.1 in \citet{van1991differentiable} gives that $\kappa(P_\lambda) \equiv \psi(\lambda)$ is not differentiable, so that $\beta_0$ cannot be estimated at root-n rate.   
\end{ex}

\subsection{Closed range and local identifiability}

For models parametrized by a Hilbert space, it is well established that, if the data tangent space $\mathcal{R}(A)$ is closed, then the local identifiability condition $\mathcal{N}(A) \subseteq \mathcal{N}(\dot\psi)$ is sufficient for $\mathcal{I} > 0$. The result extends to models parametrized by normed spaces.
\begin{prop} \label{prop:closed_range}
	If $\mathcal{R}(A)$ is closed, $\mathcal{N}(A) \subseteq \mathcal{N}(\dot\psi) \iff \dot\psi \in \mathcal{R}(A^*) \iff \mathcal{I}>0$.
\end{prop}

The fact that, when $\mathcal{R}(A)$ is closed, ``local identifiability implies positive information" can be better understood by looking at the link established in Proposition~\ref{prop:posinfo_continuous}. Recall that positive information is intimately related to continuity of $\dot\psi_\mathcal{Q}A_\mathcal{Q}^{-1}$.

When $\mathcal{R}(A)$ is closed, it turns out that $A_\mathcal{Q}^{-1}$ is continuous. This comes from Banach's Inverse Mapping Theorem \citep[Th.~1 in Section~6.4]{luenberger1997optimization}. Note that, in this case, one needs $V$ to be a Banach space: it ensures that the quotient space $\mathcal{Q}$ is also a Banach space, so that the Inverse Mapping Theorem can be applied \citep[Th.~1.7.7]{megginson1998introduction}. Then, under local identifiability, $\dot\psi_\mathcal{Q}$ is well-defined and continuous, hence so is $\dot\psi_\mathcal{Q} A_\mathcal{Q}^{-1}$. This gives positive semiparametric Fisher Information (see Proposition~\ref{prop:posinfo_continuous}). This argument is illustrative, the proof of Proposition~\ref{prop:closed_range} does not require $V$ to be a Banach space.

\section{Applications}
\label{sec:applications}

In this section, I provide the details regarding estimation of the average of a known transformation (Example~\ref{ex:mean}) and solve the $q=1$ case. I also solve the problem of estimation of the density at a point (Example~\ref{ex:density}).

\begin{ex}[continues=ex_mean]

I start by filling the gaps left in the preceding discussion. I assume that the model is well specified, that is, $fp_0 \in L_1(\mu)$ for every $f \in L_q(P_0)$. Then, the parameter space satisfies $\Lambda = \{ \lambda \in L_1(P_0) \colon \E_0[|f\lambda |] < \infty, \forall f \in L_q(P_0) \}$. To see this, note that $fp(1+\lambda) \in L_1(\mu) \iff E_0[|f\lambda|] < \infty$ follows from the triangle inequality applied to $|fp_0(1+\lambda)|$ and $|f p_0 \lambda|=|f p_0(1 + \lambda - 1)|$. The following proposition characterizes the parameter space of the model: 
\begin{prop}{\label{prop:duality}}
	Consider that $q \in [1, 2]$. Then, $\Lambda = L_{q'}(P_0)$, with $1/q + 1/q' = 1$.
\end{prop}

In this example, the tangent space is the whole parameter space, $\mathcal{T} = L_{q'}(P_0)$. The gradient of the functional $\psi(\lambda)= \int g (1+\lambda)p_0d\mu$ is $\dot\psi\alpha = \int g \alpha p_0 d\mu$. To see this, take a sequence $\lambda_t$ in $L_{q'}(P_0)$ such that $\lambda_t/t \to \alpha$ as $t \downarrow 0$. Then, by H\"older's inequality,
\begin{equation}
		\left| \frac{\psi(\lambda_t)-\psi(0)}{t} - \dot\psi\alpha \right| = \left| \int (\lambda_t/t - \alpha) g p_0 d\mu \right| \leq \norm{\lambda_t/t - \alpha}_{q'} \norm{g}_{q} \to 0,
\end{equation}
as $t \downarrow 0$, where $\norm{\cdot}_q$ denotes the $L_q(P_0)$ norm. Moreover, H\"older's inequality guarantees that the gradient $\dot\psi$ is $L_{q'}(P_0)$-continuous.

The next proposition shows mean-square differentiability of $p_\lambda \equiv p_0(1+\lambda)$, giving the shape of the score operator $A\alpha = \alpha$.
\begin{prop}\label{prop:est_mean_score}
	Let $\lambda_t$ be a sequence in $L_{q'}(P_0)$ with $q' \in [2, \infty)$ such that $\lambda_t/t \to \alpha$ as $t \downarrow 0$. Then,
	\begin{equation}
		\int \left( \frac{\sqrt{p_0(1+\lambda_t)} - \sqrt{p_0}}{t} - \frac{\alpha}{2}  \sqrt{p_0} \right)^2 d\mu \to 0, \quad \text{as } t \downarrow 0.
	\end{equation}
\end{prop}
Note that $A \colon L_{q'}(P_0) \to L_2(P_0)$ is continuous since, having $P_0(\mathcal{X})=1$, H\"older's inequality gives
\begin{equation}
	\norm{A\alpha}_2 = \norm{\alpha^2}_{1}^{1/2} \leq P_0(\mathcal{X})^{q'/(q'-2)} \norm{\alpha^ 2}_{q'/2}^{1/2} = \norm{\alpha}_{q'}.
\end{equation}

I now discuss the $q=1$ case, so that $g \in L_1(P_0)$ and $\mathcal{T} = \Lambda = L_\infty(P_0)$. I illustrate how Theorem~\ref{thm:info_range} can be applied to show that semiparametric information for $\beta_0$ is zero. The dual space of $L_\infty(P_0)$ is the space of bounded additive signed measures that are absolutely continuous w.r.t. $P_0$, which is denoted by $\operatorname{ba}(P_0)=L_\infty(P_0)^*$ \citep[Th.~IV.8.16]{dunford1988linear}. That is, for every continuous linear functional $\phi\colon L_\infty(P_0)\to \mathbb{R}$, there exists a bounded additive signed measure $\nu$ such that
	\begin{equation}{\label{eq:ba_repre}}
		\phi\alpha=\int \alpha d\nu=\int \alpha \partial\nu p_0 d\mu, \quad \forall \alpha \in L_\infty(P_0).
	\end{equation}
	In the above display, $\partial\nu$ denotes the Radon-Nikodym density of $\nu$ w.r.t. $P_0$. By comparing the above display with $\dot\psi\alpha = \int \alpha g p_0 d\mu$, it follows that the measure representing $\dot\psi$ is the one with density $g$ w.r.t. $P_0$.

		The adjoint score operator $A^*\colon L_2(P_0)  \to \operatorname{ba}(P_0)$ is characterized by the identity $\langle \alpha, A^*\delta \rangle = \langle A\alpha, \delta \rangle$ for every $\alpha \in L_\infty(P_0)$ and $\delta \in L_2(P_0)$. The Riesz–Fr\'echet Representation Theorem expresses $\langle A\alpha, \delta \rangle$ as the $L_2(P_0)$ inner product $\int A\alpha \delta p_0 d\mu$. Hence, since $A\alpha=\alpha$,
	\begin{equation}
		\langle \alpha, A^*\delta \rangle = \langle A\alpha, \delta \rangle \iff \int \alpha \partial(A^*\delta) p_0 d\mu = \int \alpha \delta p_0 d\mu,
	\end{equation}
	where $\partial(A^*\delta)$ denotes the Radon-Nikodym density of the measure representing the continuous linear functional $A^*\delta$ (c.f. equation~\eqref{eq:ba_repre}). That is, $A^*$ sends $\delta \in L_2(P_0)$ to the measure with density $\delta$ w.r.t. $P_0$.

	The conclusion from the above discussion is that $\mathcal{R}(A^*)$ consists of all the measures that have a square-integrable density w.r.t. $P_0$. Hence, if $g \in L_1(P_0) \setminus L_2(P_0)$, then $\dot\psi$ is not in $\mathcal{R}(A^*)$. Therefore, by Theorem~\ref{thm:info_range}, semiparametric information for $\beta_0$ is zero and it is not possible to estimate $\beta_0$ at root-n rate. 
\end{ex}

\begin{ex}[continues=ex:density]
	The local deviation parameter space considered for the estimation of $\beta_0 \equiv p_0(x)$ in the model $\mathcal{P} \equiv \mathcal{C}(\mathcal{X})$ is based on the following proposition:
	\begin{prop} \label{prop:local_deviations}
		Suppose that $\mathcal{X}$ is a locally compact Hausdorff space, $p_0 \in \mathcal{C}(\mathcal{X})$, and $x \in \operatorname{supp}(P_0)$. Then, there exist compact neighborhoods $C, K \subseteq \mathcal{X}$, an open set $U \subseteq \mathcal{X}$, and a continuous function $u\colon \mathcal{X} \to [0,1]$ such that
		\begin{itemize}
			\item $ x \in C \subseteq U \subseteq K \subseteq \operatorname{supp}(P_0)$ and
			\item $u(y) = 1$ if $y \in C$ and $u(y) = 0$ if  $y \in \mathcal{X} \setminus U$.
		\end{itemize}
	\end{prop}
	I consider densities of the form $p_\lambda \equiv p_0 + u\lambda$, with $\lambda \in \Lambda \equiv \mathcal{C}(K)$. Since $u(y)\lambda(y) = 0$ if $y \notin K$, one can arbitrarily extend $\lambda$ from $K$ to $\mathcal{X}$. Moreover, $u\lambda$ is a continuous function, so the set $\{p_\lambda \colon \lambda \in \Lambda \} \subseteq \mathcal{P}$. Even if the local deviations do not span the whole model, they are sufficient to account for variation in the neighborhood $C$ of $x$. Indeed, for any $p \in \mathcal{P} \equiv \mathcal{C}(\mathcal{X})$, we can take the restriction to $K$ of its deviation from $p_0$, namely $ \lambda = (p - p_0)|_K \in \mathcal{C}(K)$, which gives $p(y) = p_\lambda(y)$ for every $y \in C$, since $u(y) = 1$ if $y \in C$. 

	The tangent space of the model is $\mathcal{T} = \mathcal{C}(K)$. The functional of interest is $\psi(\lambda) = p_0(x) + u(x)\lambda(x) = p_0(x) + \lambda(x)$, since $ x \in C \Rightarrow u(x) = 1$. Its gradient is $\dot\psi\alpha = \alpha(x)$, since
	\begin{equation}
		\left| \frac{\psi(\lambda_t)-\psi(0)}{t} - \dot\psi\alpha \right| = \left| \frac{\lambda_t(x)}{t} - \alpha(x) \right| \leq \norm{\lambda_t/t -\alpha}_\infty \to 0, \quad \text{as } t \downarrow 0.
	\end{equation}

	The next proposition shows mean-square differentiability of $p_\lambda \equiv p_0+u\lambda$, giving that the score operator $A\colon \mathcal{C}(K) \to L_2(P_0)$ has shape $A\alpha = u\alpha/p_0$ and is continuous.
	\begin{prop} \label{prop:est_dens_score}
		Suppose that $\mu$ is locally finite. Let $\lambda_t$ be a sequence in $\mathcal{C}(K)$ such that $\lambda_t/t \to \alpha$ in the supremum norm as $t \downarrow 0$. Then,
		\begin{equation}
			\int \left( \frac{\sqrt{p_0 + u\lambda_t} - \sqrt{p_0}}{t} - \frac{u\alpha}{2\sqrt{p_0}}  \right)^2 d\mu \to 0, \quad \text{as } t \downarrow 0.
		\end{equation}
		Moreover, there exists a constant $0 < C < \infty$ such that $\norm{u\alpha/p_0}_2 \leq C \norm{\alpha}_\infty$.
	\end{prop}

	To characterize when is information for estimating $\beta_0 \equiv p_0(x)$ positive, one needs to study the adjoint score operator $A^* \colon L_2(P_0) \to \mathcal{C}(K)^*$. The Riesz-Markov-Kakutani Representation Theorem states that the dual space of $\mathcal{C}(K)$ is the space of finite regular Borel measures on $K$, denoted $\operatorname{rca}(K)$ \citep[Th.~2.14]{rudin1970real}. That is, for every continuous linear functional $\phi\colon \mathcal{C}(K) \to \mathbb{R}$, there exists a finite regular Borel measure $\nu$ such that
	\begin{equation}
		\phi\alpha = \int \alpha d\nu, \quad \forall \alpha \in \mathcal{C}(K).
	\end{equation}
	This means that, for every $r\in L_2(P_0)$, there is a finite regular Borel measure $\nu_r$ that represents $A^*r$ on $\operatorname{rca}(K)$. By the definition of the adjoint score operator, for every $\alpha \in \mathcal{C}(K)$ and $r \in L_2(P_0)$,
	\begin{equation}
		\langle \alpha, A^*r \rangle = \langle A\alpha, r \rangle \iff \int \alpha d\nu_r = \int A\alpha r p_0 d\mu = \int \alpha u r d\mu.
	\end{equation}
	That is, $d\nu_r = u r d\mu$. In words, $A^*$ maps $r \in L_2(P_0)$ to the measure with Randon-Nikodym density $u r$ w.r.t. $\mu$. Therefore, $\mathcal{R}(A^*)$ consists of all the measures that are absolutely continuous w.r.t. $\mu$ and have a density of the form $u r$ with $r \in L_2(P_0)$.

	Since $\dot\psi\alpha = \alpha(x)$, the representer of $\dot\psi$ on $\operatorname{rca}(K)$ is the Dirac measure $\delta_x$ at $x$, namely,
	\begin{equation}
		\delta_x(B) = \begin{cases}
			1, & \text{if } x \in B, \\
			0, & \text{if } x \notin B,
		\end{cases} \text{ for every Borel set } B.
	\end{equation}
	According to Theorem~\ref{thm:info_range}, semiparametric information for estimating $\beta_0$ is positive if and only if the Dirac measure $\delta_x$ is absolutely continuous w.r.t. $\mu$. This condition requires $\mu(\{x\}) > 0$. For instance, information is positive if the support of $P_0$ is finite, since one can take $\mu$ to be the counting measure, so that $\mu(\{x\}) = 1$. Information is also positive if there is bunching at $x$. On the other hand, if $\mu$ is the Lebesgue measure, information is zero. In this last case, Theorem~3.1 in \citet{van1991differentiable} gives that $\kappa(P_\lambda) \equiv \psi(\lambda)$ is not differentiable, so that $\beta_0$ cannot be estimated at root-n rate.
	
	I close by noting that $\mu(\{ x\}) > 0$ is enough for positive information. In that case, the Dirac measure $\delta_x$ is absolutely continuous w.r.t. $\mu$ and its density is $ur$ with 
	\begin{equation}
		r(y) = \begin{cases}
			1/\mu(\{ x\}), & \text{if } y = x, \\
			0, & \text{if } y \neq x.
		\end{cases}
	\end{equation}
	Clearly, $r$ is square-integrable w.r.t. $P_0$.
\end{ex}

\section{Conclusion}
\label{sec:conclusion}

This paper shows that, for models parametrized by arbitrary normed spaces, positive semiparametric information is equivalent to the gradient of the functional of interest being in the range of the adjoint score operator. This result extends the classical Hilbert-space result in \citet[Th.~4.1]{van1991differentiable} to general normed spaces. The paper also proposes a normed-space modelling framework for estimation of the average of a known transformation and the density at a point. The main results are applied to these problems, showing that, positive information is equivalent to the known transformation having finite variance and the density being estimated at a point with positive mass, respectively. 


\newpage
\addcontentsline{toc}{section}{References}  
\makeatletter
\makeatother
\bibliographystyle{apalike}
\bibliography{references}

\cleardoublepage
\begin{appendices}

\section{Proofs}

\begin{proofc}[Proposition~\ref{prop:posinfo_continuous}]
	
    Consider first that $\mathcal{I}>0$. For every $\delta \in \mathcal{R}(A)$, there exists an $\alpha \in \mathcal{T}$ such that $\alpha =A^{-1}\delta$. Hence, information in direction $\alpha$ can be written in terms of $\delta$ as
	\begin{equation}
		\mathcal{I}(\alpha)=\frac{\norm{\delta}^2_{2}}{|\dot{\psi} A^{-1}\delta|^2}.
	\end{equation}
    Then, for every $\delta \in \mathcal{R}(A)$, since $\mathcal{I}(\alpha) \geq \mathcal{I}$,
	\begin{equation}
		 \frac{\norm{\delta}^2_{2}}{|\dot{\psi} A^{-1}\delta|^2} \geq \mathcal{I} \Rightarrow |\dot{\psi} A^{-1}\delta| \leq \frac{1}{\sqrt{\mathcal{I}}} \norm{\delta}_{2}.
	\end{equation}
	Therefore, $\dot\psi A^{-1}$ is continuous. 
	
	For the other implication, if $\dot\psi A^{-1}$ is continuous, there exists a $0 < C<\infty$ such that $|\dot{\psi} A^{-1}\delta| \leq C \norm{\delta}_{2}$ for every $\delta\in\mathcal{R}(A)$. Then, for every $\alpha\in\mathcal{T}$, since $A\alpha \in \mathcal{R}(A)$,
	\begin{equation}
		|\dot\psi\alpha| = |(\dot\psi  A^{-1}) A \alpha| \leq C\norm{A\alpha}_{2} \Rightarrow \frac{\norm{A\alpha}_{2}^2}{|\dot\psi\alpha|^2} \geq \frac{1}{C^2}.
	\end{equation}
	Thus $\mathcal{I} \geq 1/C^2 >0$.
\end{proofc}


\begin{proofc}[Proposition~\ref{prop:range_continuous}]

	Suppose first that $\dot\psi\in\mathcal{R}(A^*)$. This means that there exits $\delta_\psi^*\in L_2(p_0)^*$ such that $A^*\delta_\psi^*=\dot\psi$. Then, for every $\delta\in\mathcal{R}(A) \subseteq L_2(p_0)$:
	\begin{equation}
		\begin{aligned}
			|\dot\psi A^{-1}\delta| &\equiv   |\langle A^{-1}\delta, \dot\psi\rangle|=|\langle A^{-1}\delta, A^*\delta_\psi^*\rangle|=|\langle AA^{-1}\delta, \delta_\psi^*\rangle| \\
			& =|\langle \delta, \delta_\psi^*\rangle| \leq \norm{\delta_\psi}_{2}\norm{\delta}_{2},
		\end{aligned}
	\end{equation}
    where $\delta_\psi \in L_2(p_0)$ is the Riesz representer of the continuous functional $\delta_\psi^*\colon L_2(p_0) \to \mathbb{R}$ \citep[see Th.~2 in][p.~109]{luenberger1997optimization}.

	Suppose now that $\dot\psi A^{-1}\colon \mathcal{R}(A)\to\mathbb{R}$ is continuous. $\dot\psi A^{-1}$ is not in $L_2(p_0)^*$, since it is not defined on the whole $L_2(p_0)$. Neverthelless, by the Hahn-Banach Extension Theorem \citep[Th.~1.9.6]{megginson1998introduction}, there exists a linear and continuous functional $\delta^*\colon L_2(p_0) \to \mathbb{R}$ that extends $\dot\psi A^{-1}$. This functional satisfies $\langle \delta, \delta^* \rangle = \langle \delta, \dot\psi A^{-1} \rangle$ for every $\delta\in\mathcal{R}(A)$.

    Take $A^*\delta^* \in \mathcal{T}^*$. For any $\alpha\in\mathcal{T}$, since $A\alpha \in \mathcal{R}(A)$, 
    \begin{equation}
		\langle \alpha, A^*\delta^* \rangle = \langle A\alpha, \delta^* \rangle = \langle A\alpha, \dot\psi A^{-1} \rangle = \langle A^{-1}A\alpha, \dot\psi \rangle =  \langle \alpha, \dot\psi \rangle.
	\end{equation}
    This shows that $(A^*\delta^*)\alpha = \dot\psi\alpha$ for every $\alpha\in\mathcal{T}$, that is, $A^*\delta^* = \dot\psi$. Therefore, $\dot\psi \in \mathcal{R}(A^*)$.
	
\end{proofc}


\begin{proofc}[Lemma~\ref{lma:properties}] \mbox{} \vspace{0.5em}
	
	\begin{enumerate}
		\item Regarding $A_\mathcal{Q}$:
			\begin{itemize}
			\item Since $\alpha'\in [\alpha]$ iff $A\alpha = A\alpha'$, it holds that $A_\mathcal{Q}[\alpha]\equiv A\alpha=A\alpha' \equiv A_\mathcal{Q}[\alpha']$. Hence, $A_\mathcal{Q}$ is well defined. 
			
			\item Linear: $A_\mathcal{Q}(a[\alpha]+b[\alpha']) \equiv A_\mathcal{Q}[a\alpha + b\alpha']=A(a\alpha + b\alpha')=aA\alpha + bA\alpha'=aA_\mathcal{Q}[\alpha] + bA_\mathcal{Q}[\alpha']$.
			
			\item Continuous: Let $\pi\colon\mathcal{T}\to\mathcal{Q}$ be the quotient map $\pi\alpha \equiv [\alpha]$. Consider an open set $U \subseteq L_2(p_0)$. It holds that $\pi(A^{-1}(U)) \subseteq A_\mathcal{Q}^{-1}(U)$. Indeed, if $[\alpha]\in \pi(A^{-1}(U))$, then $\alpha\in A^{-1}(U)$ and thus $A\alpha\in U$. This implies that $A_\mathcal{Q}[\alpha]\in U$ so $[\alpha]\in A_\mathcal{Q}^{-1}(U)$. Moreover, by continuity of $A$, $A^{-1}(U)$ is open. This, on top of Lemma~1.7.11 in \citet{megginson1998introduction}, means that, $\forall [\alpha] \in A_\mathcal{Q}^{-1}(U)$, one can find an open ball centered at $[\alpha]$ and completely contained in $A_\mathcal{Q}^{-1}(U)$.
			
			\item One-to-one: If $A_\mathcal{Q}[\alpha]=A_\mathcal{Q}[\alpha']$, then $A\alpha=A\alpha'$, so $\alpha-\alpha'\in\mathcal{N}(A)$. Therefore, $[\alpha]=[\alpha']$.
			
			\item Onto $\mathcal{R}(A)$: For $\delta\in \mathcal{R}(A)$, there exists an $\alpha\in \mathcal{T}$ such that $A\alpha=\delta$. Then, $A_\mathcal{Q}[\alpha]=\delta$.
		\end{itemize}
	
	\item Take $\alpha\in\mathcal{T}$ and $\alpha'\in [\alpha]$. This means that $\alpha' - \alpha \in \mathcal{N}(A)$. Since $\mathcal{N}(A) \subseteq \mathcal{N}(\dot\psi)$, then $\alpha' - \alpha \in \mathcal{N}(\dot\psi)$, so $\dot\psi\alpha' = \dot\psi\alpha$. Therefore, $\dot\psi_\mathcal{Q}[\alpha] \equiv \dot\psi\alpha = \dot\psi\alpha' \equiv \dot\psi_\mathcal{Q}[\alpha']$, so $\dot\psi_\mathcal{Q}$ is well defined. Linearity and continuity follow as above.
	
	\item Let
	 \begin{equation}
        \mathcal{I}_\mathcal{Q}\left([\alpha]\right) \equiv \frac{ \norm{A_\mathcal{Q}[\alpha]}_{2}^2}{\left\lvert\dot\psi_\mathcal{Q}[\alpha]\right\lvert^2}.
	\end{equation}
	Since $A_\mathcal{Q}[\alpha]\equiv A\alpha$ and $\dot\psi_\mathcal{Q}[\alpha]\equiv\dot\psi\alpha$, I have that $\mathcal{I}(\alpha)=\mathcal{I}_\mathcal{Q}([\alpha])$ for every $\alpha\in\mathcal{T}$. This implies that $\mathcal{I}=\mathcal{I}_\mathcal{Q}$. To see it, note that for every $[\alpha]\in \mathcal{Q}$, $ \mathcal{I}_\mathcal{Q}\left([\alpha]\right) = \mathcal{I}(\alpha) \geq \mathcal{I}$, since $\mathcal{I}$ is the infimum. Therefore, $\mathcal{I}$ is a lower bound for $\mathcal{I}_\mathcal{Q}([\alpha])$ and thus $\mathcal{I}_\mathcal{Q} \geq \mathcal{I}$. Likewise, for every $\alpha\in\mathcal{T}$, $\mathcal{I}(\alpha)=\mathcal{I}_\mathcal{Q}([\alpha]) \geq \mathcal{I}_\mathcal{Q}$. So $\mathcal{I}_\mathcal{Q} \leq \mathcal{I}$.
	
	\item Recall that $A_\mathcal{Q}^*\colon L_2(p_0)^*\to \mathcal{Q}^*$. I show that $\dot\psi_\mathcal{Q}\in\mathcal{R}(A_\mathcal{Q}^*) \Rightarrow \dot\psi\in\mathcal{R}(A^*)$. The proof of the other direction follows \textit{mutatis mutandis}.
	
	 If $\dot\psi_\mathcal{Q}\in\mathcal{R}(A_\mathcal{Q}^*)$, there exists a $\delta^*\in L_2(p_0)^*$ such that $\dot\psi_\mathcal{Q}=A_\mathcal{Q}^*\delta^*$. Consider $A^*\delta^* \in \mathcal{T}^*$. For every $\alpha\in\mathcal{T}$, 
	\begin{equation}
		\langle \alpha, A^*\delta^* \rangle = \langle A\alpha, \delta^* \rangle =\langle A_\mathcal{Q}[\alpha], \delta^* \rangle =\langle [\alpha], A_\mathcal{Q}^*\delta^* \rangle = \langle [\alpha], \dot\psi_\mathcal{Q} \rangle = \langle \alpha, \dot\psi \rangle.
	\end{equation}
	This means that $A^*\delta^*=\dot\psi$ and therefore $\dot\psi\in\mathcal{R}(A^*)$. 
	
	\end{enumerate} 
 
\end{proofc}


\begin{proofc}[Theorem~\ref{thm:info_range}]
	
	The implication $\dot\psi\in\mathcal{R}(A^*) \Rightarrow \mathcal{I}>0$ follows the same reasoning as in \citet{van1991differentiable}.  If $\dot\psi \in \mathcal{R}(A^*)$, then there exists an $\delta_\psi^* \in L_2(p_0)^*$ such that $\dot\psi = A^*\delta_\psi^*$. Hence, $\dot\psi\alpha\equiv \langle \alpha, \dot\psi \rangle = \langle \alpha, A^*\delta_\psi^* \rangle = \langle A\alpha, \delta_\psi^* \rangle$. Note that $\delta_\psi^*\colon L_2(p_0) \to \mathbb{R}$ is continuous with norm $\norm{\delta_\psi}_2$, being $\delta_\psi \in L_2(p_0)$ its Riesz representer \citep[see Th.~2 in][p.~109]{luenberger1997optimization}. Then, for every $\alpha\in\mathcal{T}$,
	\begin{equation}
		 |\dot\psi\alpha|=|\langle A\alpha, \delta_\psi^* \rangle| \leq \norm{\delta_\psi}_{2} \norm{A\alpha}_{2}.
	\end{equation}
	This leads to
	\begin{equation}
		\mathcal{I}(\alpha) \equiv \frac{\norm{A\alpha}_{2}^2}{|\dot\psi\alpha|^2} \geq \norm{\delta_\psi}_{2}^{-2} > 0.
	\end{equation}
	Therefore, $\mathcal{I} \geq \norm{\delta_\psi}_{2}^{-2} > 0$.
	
	Consider now that $\mathcal{I}>0$. It is easy to see that this implies the local indentifiability condition $\mathcal{N}(A)\subseteq \mathcal{N}(\dot\psi)$. When information is positive, one has that, for any $\alpha\in\mathcal{T}$, $|\dot\psi\alpha|^2 \leq \mathcal{I}^{-1}\norm{A\alpha}_{2}^2$. Thus, if $\alpha\in\mathcal{N}(A)$, then it must be that $\alpha\in\mathcal{N}(\dot\psi)$. I am then in the conditions to apply Lemma~\ref{lma:properties}.
	
	Since $A_\mathcal{Q}$ is one-to-one, applying Propositions~\ref{prop:posinfo_continuous} and \ref{prop:range_continuous} to $A_\mathcal{Q}$ and $\dot\psi_\mathcal{Q}$ gives $\mathcal{I}_\mathcal{Q}>0 \iff \dot\psi_\mathcal{Q}\in\mathcal{R}(A_\mathcal{Q}^*)$. Therefore, by Lemma~\ref{lma:properties},
	\begin{equation}
		\mathcal{I}>0 \Rightarrow \mathcal{I}_\mathcal{Q}>0 \Rightarrow \dot\psi_\mathcal{Q}\in\mathcal{R}(A_\mathcal{Q}^*) \Rightarrow \dot\psi\in\mathcal{R}(A^*).
	\end{equation}

\end{proofc}


\begin{proofc}[Proposition~\ref{prop:closed_range}]
	
	 Throughout the proof, I use the following geometry of the tuple $(\mathcal{T}, \mathcal{T}^*)$. For $S \subset \mathcal{T}$, the annihilator of $S$ is defined as $S^\perp \equiv \{\alpha^* \in \mathcal{T}^* \colon \langle \alpha, \alpha^* \rangle = 0, \forall \alpha\in S\}$. In turn, for $S^* \subset \mathcal{T}^*$, the annihilator (in $\mathcal{T}$) of $S^*$ is defined as $\prescript{\perp}{}{S^*} \equiv \{\alpha \in \mathcal{T} \colon \langle \alpha, \alpha^* \rangle = 0, \forall \alpha^*\in S^*\}$.

	From the above definitions, it immediately follows that 
	\begin{equation}
		\begin{aligned}
			\mathcal{N}(A) \subseteq \mathcal{N}(\dot\psi) \Rightarrow \mathcal{N}(\dot\psi)^\perp \subseteq \mathcal{N}(A)^\perp, \text{ and} \\
			\mathcal{N}(\dot\psi)^\perp \subseteq \mathcal{N}(A)^\perp \Rightarrow \prescript{\perp}{}{(\mathcal{N}(A)^\perp)} \subseteq \prescript{\perp}{}{(\mathcal{N}(\dot\psi)^\perp)}.
		\end{aligned}
	\end{equation}
	Since both $\mathcal{N}(A)$ and $\mathcal{N}(\dot\psi)$ are closed subspaces of $\mathcal{T}$, Proposition~1.10.15 in \citet{megginson1998introduction} gives $\mathcal{N}(A)=\prescript{\perp}{}{(\mathcal{N}(A)^\perp)}$ and $\mathcal{N}(\dot\psi)=\prescript{\perp}{}{(\mathcal{N}(\dot\psi)^\perp)}$. This shows $\mathcal{N}(A) \subseteq \mathcal{N}(\dot\psi) \iff \mathcal{N}(\dot\psi)^\perp \subseteq \mathcal{N}(A)^\perp$.

	I now show that $\mathcal{N}(A)^{\perp} = \overline{\mathcal{R}(A^*)}^{w*}$, where $\overline{\mathcal{R}(A^*)}^{w*}$ is the closure of $\mathcal{R}(A^*)$ in the weak$^*$ topology of $\mathcal{T}^*$. The result follows from Lemma~3.1.16 and Proposition~2.6.6 in \citet{megginson1998introduction}. The lemma claims that $\mathcal{N}(A)=\prescript{\perp}{}{\mathcal{R}(A^*)}$, so $\mathcal{N}(A)^\perp=(\prescript{\perp}{}{\mathcal{R}(A^*)})^\perp$. Since $\mathcal{R(A^*)}$ is a subspace of $\mathcal{T}^*$, the proposition gives $(\prescript{\perp}{}{\mathcal{R}(A^*)})^\perp=\overline{\mathcal{R}(A^*)}^{w*}$. 
	
	The next step is to show that $\mathcal{N}(\dot\psi)^\perp \subseteq  \overline{\mathcal{R}(A^*)}^{w*} \iff \dot\psi \in \overline{\mathcal{R}(A^*)}^{w*}$. First, note that $\mathcal{N}(\dot\psi) = \prescript{\perp}{}{\operatorname{span}(\{\dot\psi \})}$, where $\operatorname{span}(\{\dot\psi\}) \equiv \{a \dot\psi: a\in \mathbb{R} \}$ is the linear span of $\dot\psi$. Indeed,
	\begin{equation}
		\alpha \in \mathcal{N}(\dot\psi) \Leftrightarrow \langle \alpha, \dot\psi \rangle = 0 \Leftrightarrow \langle \alpha, a\dot\psi \rangle = 0, \forall a \in \mathbb{R} \Leftrightarrow \alpha \in \prescript{\perp}{}{\operatorname{span}(\{\dot\psi\})}.
	\end{equation}
	Thus, $\mathcal{N}(\dot\psi)^\perp= (\prescript{\perp}{}{\operatorname{span}(\{\dot\psi\})})^\perp = \overline{\operatorname{span}(\{\dot\psi\})}^{w*}$ by Proposition~2.2.6 in \citet{megginson1998introduction}. This shows that $\mathcal{N}(\dot\psi)^\perp$ is the smallest weak$^*$-closed subspace containing $\dot\psi$. So, if $\dot\psi \in \overline{\mathcal{R}(A^*)}^{w*}$, then it must be that $\mathcal{N}(\dot\psi)^\perp \subseteq  \overline{\mathcal{R}(A^*)}^{w*}$, since $\overline{\mathcal{R}(A^*)}^{w*}$ is a weak$^*$-closed subspace containing $\dot\psi$. The other implication follows from $\dot\psi \in \mathcal{N}(\dot\psi)^\perp$.

	The above discussion shows that $\mathcal{N}(A) \subseteq \mathcal{N}(\dot\psi) \iff \dot\psi \in  \overline{\mathcal{R}(A^*)}^{w*}$. To conclude, just note that $\mathcal{R}(A)$ is closed if and only if $\mathcal{R}(A^*)$ is weakly$^*$ closed \citep[Th.~3.1.21]{megginson1998introduction}. Hence, if $\mathcal{R}(A)$ is closed, $\overline{\mathcal{R}(A^*)}^{w*}=\mathcal{R}(A^*)$. The remaining equivalence follows from Theorem~\ref{thm:info_range}.
	
\end{proofc}


\begin{proofc}[Proposition~\ref{prop:duality}] 
Let $\lambda \in L_{q'}(P_0)$. Then, by H\"older's inequality, for every $f \in L_q(P_0)$, $\norm{f\lambda}_1 \leq \norm{f}_q \norm{\lambda}_{q'} < \infty$. Moreover, since $P_0$ is a finite measure, H\"older's inequality also gives $\int |\lambda| dP_0 \leq P_0(\mathcal{X})^{1/q}\norm{\lambda}_{q'} < \infty$. Therefore, $\lambda \in \Lambda$.

The other inclusion follows from standard duality results. If $\int f \lambda dP_0 < \infty$ for every $f \in L_q(P_0)$, one can construct a sequence $f_n \in L_q(P_0)$ with $\int f_n \lambda^{q'} dP_0$ approaching $\int \lambda^{q'} dP_0$ as $n \to \infty$. Then, the Monotone or Dominated Convergence Theorems will give $\norm{\lambda}_{q'} < \infty$. That is, $\lambda \in L_{q'}(P_0)$. See, for instance, the proof of Theorem~6.16 in \citet{rudin1970real}.

\end{proofc}


\begin{proofc}[Proposition~\ref{prop:est_mean_score}]
	I show that
	\begin{equation}
		\norm{\frac{\sqrt{1+\lambda_t}-1}{t} - \frac{\alpha}{2}}_2\to 0 \text{ as } t \downarrow 0,
	\end{equation}
	which implies the result. Multiplying and dividing by the conjugate, I get
	\begin{equation}
		\begin{aligned}
				\frac{\sqrt{1+\lambda_t}-1}{t} - \frac{\alpha}{2} &= \frac{\lambda_t/t}{\sqrt{1+\lambda_t}+1} - \frac{\alpha}{2} \\ &= \frac{\lambda_t/t - \alpha}{\sqrt{1+\lambda_t}+1} + \alpha\left(\frac{1}{\sqrt{1+\lambda_t}+1} - \frac{1}{2}\right) \\
				&= \left(\lambda_t/t - \alpha \right)\eta_t + \alpha\left(\eta_t - 1/2\right),
		\end{aligned}
	\end{equation}
	where $\eta_t \equiv (\sqrt{1+\lambda_t}+1)^{-1}$. Note that $0 \leq \eta_t \leq 1$ almost surely for every $t$. I deal with the two terms separately. 
	
	First,
	\begin{equation}
		\norm{(\lambda_t/t - \alpha )\eta_t}_2 \leq \norm{\lambda_t/t - \alpha}_2 \to 0, \text{ as } t \downarrow 0,
	\end{equation}
	since $L_{q'}(P_0)$ convergence implies $L_2(P_0)$ convergence for $q'\geq 2$.
	
	Second, since $\alpha^2 (\eta_t - 1/2)^2 \leq \alpha^2/4$ almost surely for every $t$, with $\alpha^2/4$ integrable since $\alpha \in L_{q'}(P_0) \subseteq L_2(P_0)$, $\alpha^2 (\eta_t - 1/2)^2 \to 0$ pointwise as $t \downarrow 0$, and the Dominated Convergence Theorem, I get $\norm{\alpha(\eta_t - 1/2)}_2 \to 0$ as $t \downarrow 0$.    
\end{proofc}


\begin{proofc}[Proposition~\ref{prop:local_deviations}]

	First, $\operatorname{supp}(P_0) = p_0^{-1}((0,\infty))$ is open since $p_0$ is continuous. Therefore, $\operatorname{supp}(P_0)$ is a neighborhood of $x$ and, since $\mathcal{X}$ is locally compact, there exists a compact neighborhood $C$ of $x$ such that $K \subseteq \operatorname{supp}(P_0)$ \citep[Def.~1.8.1]{willard2012general}. Theorem~2.7 in \citet{rudin1970real} then guarantees the existence of a compact set $K$ and an open set $U$ such that $x\in C \subseteq U \subseteq K \subseteq \operatorname{supp}(P_0)$. Finally, the existence of the continuous function $u \colon \mathcal{X} \to [0,1]$ with $u(y)=1$ if $y\in C$  and $u(y)=0$ if $y\in \mathcal{X} \setminus U$ follows from Urysohn's Lemma \citep[Lemma~2.12]{rudin1970real}. 
	
\end{proofc}


\begin{proofc}[Proposition~\ref{prop:est_dens_score}]

	First, multiplying and dividing by the conjugate, I get $\sqrt{p_0+u\lambda_t}-\sqrt{p_0} = u\lambda/(\sqrt{p_0+u\lambda_t}+\sqrt{p_0})$. Let $\eta_t \equiv (\sqrt{p_0+u\lambda_t} + \sqrt{p_0})^{-1}$. After adding and subtracting $u\eta_t$, by the triangle inequality,
	\begin{equation} \label{eq:dens_score_decomp}
		\begin{aligned}
			\norm{\frac{\sqrt{p_0 + u\lambda_t} - \sqrt{p_0}}{t} - \frac{u\alpha}{2\sqrt{p_0}}}_{2,\mu} &\leq \norm{u \eta_t \cdot \left(\frac{\lambda_t}{t} - \alpha\right)}_{2,\mu} \\
			&+ \norm{u\alpha \cdot \left(\eta_t - \frac{1}{2\sqrt{p_0}}\right)}_{2,\mu},
		\end{aligned}
	\end{equation}
	where $\norm{\cdot}_{2,\mu}$ is the $L_2(\mu)$ norm. I deal with the two terms separately.
	
	Before continuing, note that $K \subseteq \operatorname{supp}(P_0)$, hence $p_0(y) > 0$ for every $y \in K$. By the Extreme Value Theorem, since $K$ is compact, there exists an $y^*\in K$ such that $p_0(y) \geq p_0(y^*) > 0$ for every $y \in K$ \citep[c.f.][Th.~2.10]{rudin1970real}. Let $p^* \equiv p_0(y^*)$. For $y \in K$, it holds that
	\begin{equation}
		\frac{1}{\sqrt{p_0(y)+u(y)\lambda(y)}+\sqrt{p_0(y)}} \leq \frac{1}{\sqrt{p_0(y)}} \leq \frac{1}{\sqrt{p^*}}.
	\end{equation}

	For the first term in the right-hand side of \eqref{eq:dens_score_decomp}, recalling that $U \subseteq K$, I have that
	\begin{equation}
		\norm{u \eta_t \cdot \left(\frac{\lambda_t}{t} - \alpha\right)}_{2,\mu}^2 \leq \int_U \eta_t^2 \cdot \left(\lambda_t/t - \alpha\right)^2 d\mu \leq \frac{\mu(U)}{p^*} \norm{\lambda_t/t - \alpha}_\infty^2.
	\end{equation}
	Note that $\mu(U) < \mu(K) < \infty$ since $\mu$ is locally finite and $K$ is compact (just cover $K$ with the finite-measure neighborhood of every point in $K$). Therefore, the above expression goes to zero as $t \downarrow 0$.

	For the second term in the right-hand side of \eqref{eq:dens_score_decomp}, for every $y \in U \subseteq K$, I have that
	\begin{equation}
		\begin{aligned}
			\left| \eta_t - \frac{1}{2\sqrt{p_0}}  \right| &= \frac{\left| \sqrt{p_0+u\lambda_t} - \sqrt{p_0}\right|}{2\sqrt{p_0}(\sqrt{p_0+u\lambda_t}+\sqrt{p_0})} \leq \frac{\left| \sqrt{p_0+u\lambda_t} - \sqrt{p_0}\right|}{2p_0} \\
			&= \frac{u|\lambda_t|}{2p_0(\sqrt{p_0+u\lambda_t}+\sqrt{p_0})} \leq \frac{|\lambda_t|}{2p_0\sqrt{p_0}} \leq \frac{\norm{\lambda_t}_\infty}{2p^*\sqrt{p^*}}.
		\end{aligned}
	\end{equation}
	Thus,
	\begin{equation}
		\norm{u\alpha \cdot \left(\eta_t - \frac{1}{2\sqrt{p_0}}\right)}_{2,\mu}^2 = \int_U \alpha^2 \cdot \left(\eta_t - \frac{1}{2\sqrt{p_0}}\right)^2 d\mu \leq \frac{\norm{\alpha}_\infty^2\mu(U)}{4(p^*)^3} \norm{\lambda_t}_\infty^2.
	\end{equation}
	By equation~\eqref{eq:reg_path}, with $\lambda_0 = 0$, the above expression goes to zero as $t~\downarrow~0$.

	Continuity of the score operator follows from similar arguments. I have that
	\begin{equation}
		\norm{u\alpha/p_0}_2^2 = \int_U \alpha^2/p_0 d\mu \leq \frac{\mu(U)}{p^*} \norm{\alpha}_\infty^2 \Rightarrow \norm{u\alpha/p_0}_2 \leq \sqrt{\frac{\mu(U)}{p^*}} \norm{\alpha}_\infty.
	\end{equation}
\end{proofc}

\end{appendices}

\end{document}